\documentclass[12pt,oneside,english]{amsart}
\usepackage[latin1]{inputenc}
\usepackage{amssymb}

\makeatletter
\newtheorem{theorem}{Theorem}

\newtheorem{example}{Example}
\newtheorem{definition}{Definition}

\newtheorem{corollary}{Corollary}

\usepackage{babel}

\makeatother
\begin{document}
\baselineskip=17pt

\title{Process of "Primoverization" of Numbers of the Form $a^n-1$}

\author{Vladimir Shevelev}
\address{Department of Mathematics \\Ben-Gurion University of the
 Negev\\Beer-Sheva 84105, Israel. e-mail:shevelev@bgu.ac.il}

\begin{abstract}
 We call an integer $N>1$ primover to base $a$ if it either prime or overpseudoprime to base $a$. We prove, in particular, that every Fermat number is primover to base 2. We also indicate a simple process of receiving of primover divisors of numbers of the form $a^n-1.$

\end{abstract}

\maketitle

\section{Introduction}
   Let $a>1$ be an integer. For an odd $n >1$, consider the number $r=r_a(n)$ of distinct cyclotomic cosets of $a>1$ modulo $n$  [1, pp.104-105].
Note that, if $C_1,\ldots, C_r$ are all different cyclotomic cosets of $a\mod n$, then

\begin{equation}\label{2}
\bigcup^r_{j=1}C_j=\{1,2,\ldots, n-1\},\qquad C_{j_1}\cap C_{j_2}=\varnothing , \;\; j_1\neq j_2.
\end{equation}

For the least common  multiple of $|C_1|, \ldots, |C_r|$ we have

\begin{equation}\label{3}
[|C_1|,\ldots,|C_r|]= h,
\end {equation}
where $h=h_a(n)$ is the multiplicative order of $a$ modulo $n.$
(This follows easily, e.g., from Exercise 3, p. 104 in \cite{2}).

It is easy to see that for odd prime $p$ we have

\begin{equation}\label{4}
|C_1|=\ldots=|C_r|
\end {equation}

such that

\begin{equation}\label{5}
p= rh + 1.
\end {equation}

\begin{definition}
We call odd composite number $n$ overpseudoprime to base $a$ $(n\in\mathbb{S}_a)$ if
\newpage
\begin{equation}\label{6}
n=r_a(n) h_a(n)+1.
\end {equation}

\end{definition}

Let $n$ be odd composite number with the prime factorization
\begin{equation}\label{7}
n=p_1^{l_1}\cdots p_k^{l_k}.
\end {equation}
In [3] we proved the following criterion.
\begin{theorem}\label{t1}
The composite number $n$ (6) is overpseudoprime to base $a$ if and only if for all nonzero vectors $(i_1, \ldots, i_k)\leq (l_1, \ldots, l_k)$ we have
\begin{equation}\label{8}
h_a(n)=h_a(p_1^{i_1}\cdots p_k^{i_k}).
\end {equation}
\end{theorem}
\begin{corollary}\label{1} Every two overpseudoprimes $n_1$ and $n_2$ for which $h_a(n_1)\neq h_a(n_2)$ are
coprimes.
 \end{corollary}\label{1}
Notice that, every overpseudoprime to base $a$ is always  superpseudoprime  and
strong pseudoprime to the same base  ( see Theorems 12-13 in [3]).
Besides,in [3] we proved, in particular, the following result.
\begin{theorem}\label{2} If $p$ is prime then  $2^p-1$ is  either prime or overpseudoprime to base 2.

 \end{theorem}\label{2}
Note that, prime divisors of overpseudoprime $n$  are primitive divisors of $2^{h(n)}-1.$Thus, all overpseudoprimes to base 2 are generated by the set of primitive divisors of the sequence $(2^n-1)_{n\geq1}.$\newline
   Denote by $Ov_2(x)$ the number of overpseudoprimes to base 2 not exceeding $x.$ In the next paper of our cycle  we shall   prove that
\begin{equation}\label{9} Ov_2(x)=o(x^\varepsilon),
\end {equation}
where $\varepsilon>0$ is arbitrary small for sufficiently large $x.$
\begin{definition}An integer which is either prime or overpseudoprime to base $a$ we call primover to the same base.
\end{definition}
The following theorem complements Theorem 2.
\begin{theorem}\label{3} If $n$ is primover to base 2, then $2^n-1$ is primover to base 2 if and only if $n$ is prime.
\end{theorem}\label{3}
\bfseries Proof.\enskip\mdseries If $n$ is prime then, by Theorem 2, $2^n-1$ is primover to base 2. Let $n$ be not prime.  If it contains two different prime divisors $p<q,$ then $h_2(2^p-1)=p<q=h_2(2^q-1).$ This means that $2^n-1$, which is multiple of $2^p-1$ and of $2^q-1$, is neither prime, nor overpseudoprime to base 2. At last, if
$n=p^k$ then $h_2(2^p-1)=p<h_2(2^{p^2}-1)$ and again $n$ is not primover \newpage to base 2. $\blacksquare$ \newline
     In this paper we indicate a very simple process of a constructive construction of (large) primovers to base $a$ from numbers of the form $a^n-1.$
\section{Process of "primoverization"}
  We start with the following result.
  \begin{theorem}\label{4} For $n\geq1$ every Fermat number $ F_n=2^{2^{n-1}}+1$ is primover to base 2.
\end{theorem}
\bfseries Proof.\enskip\mdseries  Let $d>1$ be a divisor of $F_n.$ Then $2^{2^{n-1}}\equiv -1\pmod d$ and
$2^{2^{n}}\equiv 1\pmod d.$ Thus, $h_2(d)\leq 2^d.$  But $h_2(d)|2^n.$ Therefore, $h_2(n)=2^i, i\leq n.$
If $i\leq n-1$ then d divides $2^{2^{i}-1}$ and does not divide $2^{2^{i-1}-1}.$ Hence, $d|2^{2^{i-1}}+1=F_i$. It is impossible since for $i\leq n-1$ we have $(F_i, F_n)=1.$ Therefore, for every
divisor $d>1$ of $F_n$ we have $h_2(d)=2^n$ and, by  Theorem 1, if $F_n$ is not prime, then it is overpseudoprime to base 2. $\blacksquare$
    Thus, if the well known Hardy-Wright conjecture about the finiteness of the number of Fermat primes
is true, then the expression $2^{2^{n-1}}+1$ gives since some $n$ only overpseudoprimes to base 2. On the other hand , by Theorem 3, the expression $2^{F_n}-1$ gives since the corresponding $n$ neither primes nor overpseudoprimes to base 2. It is interesting also to notice that , if $Pr_2(x)$ denotes the number of primovers to base 2 not exceeding $x$ , then by (8) we have $Ov_2(x)=\pi(x)+o(x^\varepsilon).$ In spite
of this very small difference, the problem of finding of a useful formula which gives only primes is unsolved until now, while the simple and very useful in mathematics Fermat construction gives only primovers to base 2.
\begin{example}For $n=6$ we obtain the following primover to base 2: 4294967297. It is 2315-th strong pseudopime to base 2 (see [4, Seq. A001262]).
\end{example}
Quite analogously one can prove the following generalization of Theorem 4.
\begin{theorem}\label{5}For $n\geq1$ and even $a\geq2$ every generalized Fermat number $ F_n^{(a)}=a^{2^{n-1}}+1$ is primover to base $a.$
\end{theorem}
    The equality $a^{2^{n-1}}+1=\frac {a^{2^{n}}-1} {a^{2^{n-1}}-1}$ promts a way of "primoverization" of
numbers of the form $2^m-1$ and,  in more general, of the form $a^m-1.$\newpage
\begin{theorem}\label{6}If $p<q$ are primes , then the number
$$ N=\frac {(a-1)(a^{pq}-1)} {(a^p-1)(a^q-1)} $$
is primover to base $a$  if and only if  $(N,(a^p-1)(a^q-1))=1. $
\end{theorem}
 \bfseries Proof.\enskip\mdseries  Notice, that $ (a^{p}-1, a^{q}-1)=a-1$ and , therefore, $N$ is integer. Let  $(N,(a^p-1)(a^q-1))=1 $ and $d>1$ be a divisor of $N.$ Then $d$ is a divisor of $a^{pq}-1$ and,
 hence, $h_a(d)|pq.$ But if $h_a(d)=p$ then $d|a^p-1.$ This contradicts to the condition.
  By the same arguments $h_a(d)\neq q$. Thus, $h_a(d)=pq$ and does not depend on $d.$  Therefore, using Theorem 1, we conclude that $N$ is primover. Conversly, let $N$ be primover.Then for every its divisor $d>1$ we have $h_a(d)=pq.$
  Therefore, $d$ divide neither $a^{p}-1$ nor $a^{q}-1.$ Thus, $(N,(a^p-1)(a^q-1))=1.\blacksquare$
 \begin{example}For $a=2, p=5, q=7$  we obtain a primover 8727391 which is the 150-th strong pseudopime to base 2 (see [4, Seq. A001262]).
 \end{example}
 Quite analogously one can prove the following result.
\begin{theorem}\label{7}If $p$ is prime, then
$$ N=\frac {a^{p^n}-1} {a^{p^{n-1}}-1} $$ is primover to base $a$ if and only if  $(N, a^{p^{n-1}}-1)=1.$
\end{theorem}
\begin{example}For $a=2,p=5,n=2$ we obtain a primover 1082401 which is the 50-th strong pseudopime to base 2 (see [4, Seq. A001262]).
\end{example}
 Let n is a composite squarefree number and its prime factorization is $n=p_1\cdot\ldots\cdot p_t.$ Let to each combination of  $k\geq0$ elements $i_1,\cdots,i_k$ from set $(1,\ldots,t)$ corresponds the number
 $$a^{p_{i_1}\cdot\ldots\cdot p_{i_k}}-1=m_n(i_1,\cdots,i_k)$$
 (in case of $k=0$ we put $m_n=a-1$). \newline
 The following result, which could be proved using induction, generalizes Theorem 6.
 \begin{theorem}\label{8} The number
  $$N=\frac {\prod_{k\geq 0, k\equiv t \pmod 2,\enskip1\leq i_1<\ldots<i_k\leq t}m_n(i_1,\cdots,i_k) } {\prod_{k\geq 0, k\equiv t-1 \pmod 2,\enskip1\leq i_1<\ldots<i_k\leq t}m_n(i_1,\cdots,i_k)}$$
   is primover to base $a$ if and only if the number $a^{p_1\cdot\ldots\cdot p_t}-1$ is coprime to $\frac { a^{p_1\cdot\ldots\cdot p_t}-1}{N}.$
 \end{theorem}\newpage
 \begin{example} For $a=2, t=3, p_1=2, p_2=5, P_3=7$   we obtain a primover 24214051 which is the 254-th strong pseudopime to base 2 (see [4, Seq. A001262]).
 \end{example}
      Furthermore, using the construction of Theorem 7, we obtain the following generalization of Theorem 6.
\begin{theorem}\label{9}If $p<q$ are primes then, for $a\geq 1, b\geq 1,$   the number
$$ N=\frac {(a^{p^{a}q^{b}}-1)(a^{p^{a-1}q^{b-1}}-1)} {(a^{p^{a-1}q^{b}}-1)(a^{p^{a}q^{b-1}}-1)} $$
is primover to base $a$ if and only if the number $a^{p^{a}q^{b}}-1$ is coprime to
$$\frac{(a^{p^{a-1}q^{b}}-1)(a^{p^{a}q^{b-1}}-1)}{(a^{p^{a-1}q^{b-1}}-1)}.$$
\end{theorem}
      In order to write a generalization we need to represent a nonnegative integer $m<2^{k}$ in the $k$-digit
binary expansion ( probably, with some 0's before the first 1). It is known, that $m$ is called odious (evil) if the number of 1's in its binary expansion is odd (even).
      In the general case of composite number $n=p_1^{l_1}\cdot\cdots\cdot p_k^{l_k}$ it could be obtained the following result.
\begin{theorem}\label{10}  The number
\begin{equation}\label{10}
 N=\frac {\prod(a^{p_1^{l_1-i_1}\cdot\cdots\cdot p_k^{l_k-i_k}}-1)} {\prod(a^{p_1^{l_1-j_1}\cdot\cdots\cdot p_k^{l_k-j_k}}-1)},
\end{equation}
where $(i_1,\ldots, i_k)$ runs the vectors of binary digits of all k-digit evil numbers from $[ 0, 2^{k}-1],$ while
$(j_1,\ldots, j_k)$ runs the vectors of binary digits of all k-digit odious numbers from the same segment,
is primover to base $a$ if and only if   $a^{p_1^{l_1}\cdot\cdots\cdot p_k^{l_k}}-1$ is coprime to $\frac {a^{p_1^{l_1}\cdot\cdots\cdot p_k^{l_k}}-1} {N} .$
\end{theorem}
\begin{definition}Let $n=p_1^{l_1}\cdot\cdots\cdot p_k^{l_k}.$ We call the number $N,$ which is defined by (9),the primitive primover cofactor of $a^{n}-1,$ writing $N=Pr(a^{n}-1)$.
\end{definition}

     Using induction, we obtain the following result.
\begin{theorem}\label{11}If in the expression for primover cofactor $N$ (9) to drop in all factors $(-1)'s$, then we
obtain $a^{\varphi(n)},$ where $\varphi(n)$ is the Euler function.
\end{theorem}
Finally, we find a lower estimate of the primitive primover cofactor of $a^{n}-1.$
\begin{theorem}\label{12}There exist a positive constant $C_1$ such that
$$ Pr(a^{n}-1)\geq (a^{n}-1)^{C_1/\ln\ln n} .$$
\end{theorem}\newpage
\bfseries Proof \enskip\mdseries easily follows from Theorem 11 and the well-known Landau inequality $\varphi(n)\geq Cn/\ln\ln n.\blacksquare$

\end{document}